\documentclass[twoside]{article}
\usepackage{latexsym}
\usepackage{amsfonts}
\usepackage{graphicx}

\newtheorem{theorem}{Theorem}[section]
\newtheorem{lemma}{Lemma}[section]
\newtheorem{cor}{Corollary}[section]

\newtheorem{prop}{Proposition}[section]
\newcommand{\qed}{\hfill$\Box$\par\medskip}

\newenvironment{Proof}{\noindent{\sc Proof.}}{\qed}

\def\bhag#1{\noindent
\setcounter{equation}{0}
\section{#1}
}
\def\bfgk#1{{{#1}\kern-5.5pt{#1}}}

\def\seqb{{\mathsf b}}

\def\RR{{\mathbb R}}
\def\CC{{\mathbb C}}
\def\ZZ{{\mathbb Z}}

\def\PPI{{{\rm I}\kern-1pt\Pi}}

\def\a{\alpha}
\def\b{\beta}

\def\z{{\bf z}}

\def\O{{\cal O}}

\def\C{{\mathcal C}}

\def\Z{{\cal Z}}

\def\derf#1#2{{#1}^{(#2)}}

\def\esssup{\mathop{\hbox{{\rm ess sup}}}}
\def\supp{{\mbox{{\rm supp\ }}}}

\def\be{\begin{equation}}
\def\ee{\end{equation}}
\def\bea{\begin{eqnarray}}
\def\eea{\end{eqnarray}}
\def\eref#1{(\ref{#1})}
\def\disp{\displaystyle}

\def\jsup#1#2#3{{{#1}^{(#2,#3)}}}

\def\donchitre#1#2{\vskip 6.5cm\noindent
\parbox[t]{1in}{\special{eps:#1.eps x=6.5cm y=5.5cm}}
\hbox to 7cm{}\parbox[t]{0.0cm}{\special{eps:#2.eps x=6.5cm y=5.5cm}}}

\title{Polynomial operators and local smoothness classes on the unit interval, II}
\author{
H.~N.~Mhaskar\thanks{The research of this author was supported, in part, by grant DMS-0605209 from the National Science Foundation and  grant W911NF-04-1-0339 from the U.S. Army Research Office.}\\
Department of Mathematics, California
State University\\
Los Angeles, California, 90032, U.S.A.\\
hmhaska@calstatela.edu}

\date{}

\begin{document}
\maketitle
\hspace{1.0cm}

\begin{abstract}
We prove the existence of quadrature formulas exact for integrating high degree polynomials with respect to Jacobi weights based on scattered data on the unit interval. We also obtain a characterization of local Besov spaces using the coefficients of a tight frame expansion.
\end{abstract}

\medskip\noindent
\bhag{Introduction}

It is well known that a major drawback of polynomial approximation is that polynomials cannot be localized; a polynomial of degree $n$ is completely determined by its values at $n+1$ points on an interval, howsoever small. Another example is the following. If $P_n^*$ is a best polynomial approximation of degree at most $n$ to the function $f(x)=|x|$ on $[-1,1]$, then there are at least $n+2$ points $y_j$ on  $[-1,1]$, where $|f(y_j)-P_n^*(y_j)|\ge cn^{-1}$, where $c$ is a positive constant independent of $n$ and the points $y_j$. Moreover, as $n\to\infty$, these points become dense on $[-1,1]$. Thus, even though the function is piecewise analytic, the only singularity of the function, namely, $x=0$, affects the quality of best approximation on the entire interval. 

One of the interesting problems in the theory of polynomial approximation is thus to construct localized polynomial approximations. There are two flavors to this problem. In the \emph{local information problem}, one has information about the target function only on a part of the interval $[-1,1]$, and wants to construct the approximation also on a neighborhood of this part, keeping the growth of the approximation in control on the rest of the interval. In the \emph{local adaptivity problem},  the data is available on the whole interval; indeed, it might even be in the form of the globally defined Fourier coefficients $\hat f(\mu;k)$ (defined in \eref{fourcoeffdef}) of the target function $f$  for a suitable measure $\mu$. Nevertheless, we wish to find a single polynomial that adapts its behavior on different parts of the interval according to the smoothness of the target function on these parts. 
This is a different problem from piecewise polynomial (or spline) approximation and adaptive approximation. Unlike in spline approximation, we want to obtain a single polynomial, which yields the benefits of the superior degree of approximation provided by polynomials on intervals where the function is smooth or even analytic. Unlike in adaptive approximation, we   do not wish first to find the location of the singularities of the target function; the determination of their locations should be a byproduct of successive polynomial approximations. The construction of these polynomials may depend upon finitely many of the  Fourier coefficients, or on the values of $f$ at certain points on the interval. In most modern applications, one does not have a control on where to choose these points. Such problems are called \emph{scattered data} problems. Thus, one of the important problems in this theory is to develop quadrature formulas similar to the Gauss--Jacobi quadarture formula, except that one should not require that the values of $f$ be at the zeros of some orthogonal polynomial, but may be at any set of points instead. In return, we cannot expect the formula based on $m$ points to be exact for polynomials of degree $2m-1$ but only for polynomials of degree proportional to $m$. 

To summarize, we wish to construct a wavelet--like expansion $\sum_k\sum_{j}a_{j,k}(f)\psi_{j,k}$ with the following properties: (1) For $k=0,1,\cdots$, $\psi_{j,k}$ is a polynomial of degree at most $2^k$, independent of $f$. (2) Each coefficient $a_{j,k}$ should be a finite linear combination of the Fourier coefficients of the target function $f$, or the values of $f$ at $\O(2^k)$ arbitratily chosen points on the interval. (3) The difference  $\left|\sum_{k=0}^n \sum_{j}a_{j,k}(f)\psi_{j,k}(x)-f(x)\right|$ should be $\O(E_{\mu;c2^n,\infty}(f))$ uniformly on $[-1,1]$ (see \eref{degapproxdef} for definition, $\mu$ being a suitable measure), and yet, near every point $x\in [-1,1]$, it should be small according to the behavior of $f$ near $x$.  (4) It should be possible to characterize the smoothness of $f$ near different points by the absolute values $|a_{j,n}(f)|$. (5) The following frame inequalities should be satisfied:
\be\label{framebds}
\sum_k\sum_{j}|a_{j,k}(f)|^2 \sim \|f\|_{\mu;2}.
\ee
The system $\{\psi_{j,k}\}$ will be called a \emph{frame}. In the case when equality holds in \eref{framebds} rather than $\sim$,  the system $\{\psi_{j,k}\}$ will be called a \emph{tight frame}.

This paper is a continuation of our paper \cite{locjacobi}, where we constructed certain localized kernels, and characterized local Besov spaces on the unit interval in terms of polynomial frames based on Jacobi polynomials. 
The aim of this paper is to fill in certain gaps left over in \cite{locjacobi}, which have come to our attention after the publication of that paper. 
In Section~\ref{backsect}, we introduce the necessary notations. In Section~\ref{framesect}, we develop the localized tight frames in the context of general mass distributions on $[-1,1]$. Essential ingredients in this theory are localized polynomial kernels and quadrature formulas.  In Section~\ref{quadsect}, we develop quadrature formulas exact for integrating high degree polynomials with respect to the Jacobi weight, based on scattered data on the interval. In Section~\ref{kernsect}, we prove localization estimates on certain kernels based on Jacobi polynomials that are more elegant than those presented in \cite{locjacobi}. Certain computational issues are addressed in Section~\ref{computesect}.

\bhag{Notations and background}\label{backsect}

For $x>0$, let $\Pi_x$ denote the class of all (algebraic) polynomials of degree at most $x$. Denoting by $\lfloor x\rfloor$ the integer part of $x$, this is the same as the class $\Pi_{\lfloor x\rfloor}$. Extending the notation in this way allows us to use a less cumbersome notation, as well as relieves us of stating certain conditions on the degrees in almost every statement. If $\mu$ is a finite positive or signed Borel measure on $[-1,1]$, the total variation measure of $\mu$ will be denoted by $|\mu|$. A point $x\in [-1,1]$ is called a point of increase of $\mu$ if $|\mu|(I)>0$ for every open interval $I$ containing $x$. The set of all points of increase of $\mu$ is called the support of $\mu$, denoted by $\supp(\mu)$. It is clear that $\supp(\mu)$ is always a closed set. The measure $\mu$ is called a mass distribution if it has infinitely many points of increase and $|\mu|([-1,1])<\infty$. If $A\subset [-1,1]$ is a $\mu$--measurable subset, and $f : A\to \CC$, we  define
$$
\|f\|_{\mu;p,A}:=\left\{\begin{array}{ll}
\disp \left\{\int_A |f(t)|^pd|\mu|(t)\right\}^{1/p}, & \mbox{ if $1\le p <\infty$,}\\
|\mu|-\esssup_{t\in A}|f(t)|, &\mbox{ if $p=\infty$.}
\end{array}\right.
$$
The class $L^p(\mu,A)$ consists of all $\mu$--measurable functions $f : A\to\CC$ for which $\|f\|_{\mu;p,A}<\infty$, with the convention that two functions are considered equal if they are equal $|\mu|$--almost everywhere.  The space $C(A)$ denotes the class of all uniformly continuous, bounded functions on $A$.   The symbol $X^p(\mu,A)$ denotes $L^p(\mu,A)$ if $1\le p<\infty$, and $C(A)$ if $p=\infty$. 
For $f\in L^p(\mu,A)$ and $x\ge 0$, we define the degree of approximation of $f$ from $\Pi_x$ by
\be\label{degapproxdef}
E_{\mu;x,p,A}(f):=\inf_{P\in\Pi_x}\|f-P\|_{\mu;p,A}.
\ee
In the sequel, if   $A=[-1,1]$, we will often omit its mention from the notations, when not required by considerations of clarity. Similarly, we will often omit the mention of the measure $\mu$ if it is the Lebesgue measure.

It is convenient to define the Besov spaces using the degrees of approximation. Characterizations in terms of the Ditzian--Totik moduli of smoothness are given in \cite{devlorbk} for some measures $\mu$. Let $0<\rho\le \infty$, $\gamma>0$, and ${\bf a}=\{a_n\}_{n=0}^\infty$ be a sequence of real numbers. We define
\be\label{seqbdef}
\|{\bf a}\|_{\rho,\gamma} := \cases{\disp \{\sum_{n=0}^\infty 2^{n\gamma\rho}|a_n|^\rho\}^{1/\rho}, & if $0<\rho<\infty$,\cr
\disp\sup_{n\ge 0} 2^{n\gamma}|a_n|, & if $\rho=\infty$.\cr}
\ee
The space of sequences ${\bf a}$ for which $\|{\bf a}\|_{\rho,\gamma}<\infty$ will be denoted by $\seqb_{\rho,\gamma}$. 
For $1\le p\le\infty$, $0<\rho\le \infty$, $\gamma>0$, the Besov space $B_{\mu;p,\rho,\gamma}:=B_{\mu;p,\rho,\gamma}$ consists of functions $f\in X^p(\mu)$ for which the sequence $\{E_{\mu;2^n,p}\}\in \seqb_{\rho,\gamma}$. We define
\be\label{besovnorm}
\|f\|_{B_{\mu;p,\rho,\gamma}}=\|f\|_{\mu;p}+\|\{E_{\mu;2^n,p}\}\|_{\rho,\gamma}.
\ee 
Let $x_0\in [-1,1]$. The local Besov space $B_{\mu;p,\rho,\gamma}(x_0)$ consists of all $f\in X^p(\mu)$ with the property that for every $C^\infty$ function $\phi$ supported on an interval containing $x_0$, $f\phi\in B_{\mu;p,\rho,\gamma}$. 
Let $\mu$ be a finite positive measure on $[-1,1]$. If $\mu$ is a mass distribution, one can use the Gram--Schmidt orthogonalization process to obtain a unique sequence of orthonormalized polynomials $p_k(\mu;x)=\gamma_k(\mu)x^k+\cdots\in \Pi_k$, $\gamma_k(\mu)>0$, $k=0,1,\cdots$, such that
$$
\int_{-1}^1 p_k(\mu;x)p_j(\mu;x)d\mu(x)= 0, \qquad \mbox{if $k\not=j$, $k, j =0,1,\cdots$},
$$
and $\int_{-1}^1p_k^2(\mu;x)d\mu(x)=1$, $k=0,1,\cdots$. It is customary to define $p_k(\mu;x)=0$ if $k<0$. If $f\in L^1(\mu)$, we may define
\be\label{fourcoeffdef}
\hat f(\mu;k)=\int_{-1}^1 f(t)p_k(\mu;t)d\mu(t), \qquad k=0,1,\cdots,
\ee
and for $m=1,2,\cdots$,
\be\label{foursumdef}
s_m(\mu;f,x):=\sum_{k=0}^{m-1}\hat f(\mu;k)p_k(\mu;x)=\int_{-1}^1 f(t)K_m(\mu;x,t)d\mu(t),
\ee
where the Christoffel--Darboux kernel $K_m(\mu;x,t):=\sum_{k=0}^{m-1}p_k(\mu;x)p_k(\mu;t)$.   In  particular,
\be\label{darbouxreprod}
\int_{-1}^1 K_m(\mu;x,t)^2d\mu(t)=K_m(\mu;x,x), \qquad x\in \RR, \ m=1,2,\cdots.
\ee
We find it convenient to introduce the notation
$$
\lambda_{m}(\mu;x)^{-1}:= K_m(\mu;x,x), \qquad x\in\RR, \ m=1,2,\cdots.
$$

 In the sequel, the symbols $c, c_1,\cdots$ will denote generic positive constants independent of the degrees of polynomials involved and the target function $f$, but may depend upon such fixed quantities in the discussion as $\mu$ and the smoothness parameters $\gamma$, $\rho$, etc. The symbol $A_1\sim A_2$ means that $c_1A_1\le A_2\le c_2A_1$.

\bhag{Polynomial frames}\label{framesect}
 We will assume in the sequel that $\mu$ is a fixed mass distribution. 

Let $S\ge 2$ be an integer, and $h^\circ :\RR\to [0,\infty)$ be a compactly supported function that can be expressed as an $S$ times iterated integral of a function of bounded variation on $\RR$. In view of the Poisson summation formula \cite{butzer}, 
$$
\sum_{k\in\ZZ}h^\circ(k/\lambda)e^{ikx} =  \lambda \sum_{k\in\ZZ}
    \int_{-\infty}^\infty 
        h^\circ(t)\exp\biggl(i\lambda (x+2k\pi)t\biggr)dt, \qquad x\in [-\pi,\pi],\ \lambda>0.
$$
It is not difficult to derive from this formula by a repeated integration by parts that (cf. \cite{trigwave})
$$
\left| \sum_{k\in\ZZ}h^\circ(k/\lambda)e^{ikx}\right| \le cV(\derf{h^\circ}{S-1})\lambda\min\{1,(\lambda|x|)^{-S}\}, \qquad x\in [-\pi,\pi], \lambda>0,
$$
where $V(\derf{h^\circ}{S-1})$ is the total variation of $\derf{h^\circ}{S-1}$ on $\RR$. The following definitions are motivated by this observation.

If $h :[0,\infty)\to \RR$ is compactly supported, let
\be\label{phikerndef}
\Phi_\lambda(\mu;h,x,t) := \sum_{k=0}^\infty h\left(\frac{k}{\lambda}\right)p_k(\mu;x)p_k(\mu;t), \qquad x,t\in\RR,\ \lambda \ge 1.
\ee
 If $\lambda <1$, then we set $\Phi_\lambda(\mu;h,x,t)=0$. It may seem strange to define it so. However, in most applications, $h(t)=0$ if $t\ge 1$. Therefore, even if $\Phi_\lambda(\mu;h,x,t)$ were defined by \eref{phikerndef} for all $\lambda>0$, we would have $\Phi_\lambda(\mu;h,x,t)=h(0)p_0(\mu)^2=\Phi_1(\mu;h,x,t)$ for all $\lambda \le 1$.  The statements of our results later (e.g., Corollary~\ref{framecor}) will be simplified with the definition as we have made it. We note that if we choose $h(t)=1$ for $0\le t<1$ and $h(t)=0$ if $t\ge 1$, then for integer $n\ge 1$, the kernel $\Phi_n(\mu;h)$ reduces to the Christoffel--Darboux kernel. In the sequel, we will assume that $\lambda\ge 0$.

 Let $S\ge 1$ be an integer. A function $h: [0,\infty)\to [0,\infty)$ will be called a \emph{multiplier mask of order $S$} if each of the following conditions \eref{multmaskcond}, \eref{localizationcond} is satisfied:
\be\label{multmaskcond}
\sup_{\lambda\ge 0,\ x\in [-1,1]}\int_{-1}^1 |\Phi_\lambda(\mu;h,x,t)|d\mu(t) <\infty,
\ee
\be\label{localizationcond}
\sup_{x,t\in [-1,1],\ |x-t|\ge\delta} |\Phi_\lambda(\mu;h,x,t)|\le c(\delta)\lambda^{-S}, \qquad \lambda, \delta>0.
\ee
If $\nu$ is a (possibly signed) measure, we define the \emph{summability operators} by
\be\label{gensigmaopdef}
\sigma_\lambda(\nu,\mu;h,f,x):=\int_{-1}^1 \Phi_\lambda(\mu;h,x,t)f(t)d\nu(t), \qquad \lambda\ge 0,\ x\in\RR.
\ee
If $\nu=\mu$, then we write $\sigma_\lambda^*(\mu;h,f):=\sigma_\lambda(\nu,\mu;h,f)$. We note that if $\lambda\ge 1$ then 
\be\label{sigmacontexp}
\sigma_\lambda^*(\mu;h,f,x)=\sum_{k=0}^\infty h\left(\frac{k}{\lambda}\right)\hat f(\mu;k)p_k(\mu;x), \qquad x\in\RR.
\ee
Since $h$ is compactly supported, $\sigma_\lambda^*(\mu;h,f,x)$ can be computed using $\O(\lambda)$ Fourier coefficients of $f$. If $\nu$ is a discrete measure that associates the mass $w_z$ with each point $z$ in a finite set $\C$, then
$$
\sigma_\lambda(\nu,\mu;h,f,x)=\sum_{\z\in\C} w_z\Phi_\lambda(\mu;h,x,z)f(z),
$$
which can be computed using the information $\{f(z)\}_{z\in\C}$. Our theory usually does not depend upon the choice of the set $\C$ and the specific weights $w_z$. Therefore, we will often use the Stieltjes notation to denote a sum of the form $\sum_{z\in\C}w_zf(z)$. Thus, we define a measure $\nu$ to associate the mass $w_z$ with each $z\in\C$, and write $\int fd\nu$ or a variant thereof to denote the finite sum.   If $1\le p\le\infty$, then we will say that $\nu\preceq_p \mu$ if $f\in X^p(\mu)$ implies that $f\in X^p(\nu)$, and $\|f\|_{\nu;p}\le c\|f\|_{\mu;p}$.

 A (possibly signed) measure $\nu$ will be called a M--Z (Marcinkiewicz--Zygmund) quadrature measure of order $n$ (for $\mu$)  if  each of the following conditions \eref{genquad}, \eref{genmzineq} is satisfied:
\be\label{genquad}
\int_{-1}^1 Pd\nu =  \int_{-1}^1 Pd\mu, \qquad P\in\Pi_n,
\ee
\be\label{genmzineq}
\|P\|_{\nu;1} \le c \|P\|_{\mu;1}, \qquad P\in\Pi_n.
\ee
(As usual, we tacitly think of the measure $\nu$ to be a member of a sequence of measures.) 

 In the rest of this section, let $h: [0,\infty)\to [0,\infty)$ be a nonincreasing function, such that $h(t)=1$ if $0\le t \le 1/2$ and $h(t)=0$ if $t\ge 1$. We will assume that $h$ is a multiplier mask of order $S$. The constants $c, c_1,\cdots$ may depend upon $h$.

The following proposition lists some immediate properties of the summability operators.
\begin{prop}\label{sigmaprop}
Let $m\ge 0$ be an integer, $\nu$ be an M--Z quadrature measure of order $3m/2-1$. Let $1\le p\le\infty$ and $f\in L^p(\mu)$.\\
{\rm (a)} For $P\in \Pi_{m/2}$, $\sigma_m(\nu,\mu; h,P)=P$.\\
{\rm (b)}  We have 
\be\label{uniformbd}
\|\sigma_m(\nu,\mu; h,f)\|_{\mu;p}\le c\|f\|_{\nu;p}, \qquad m=0,1,\cdots.
\ee
Consequently, if $\nu\preceq_p\mu$ and $f\in X^p(\mu)$, then $\|\sigma_m(\nu,\mu; h,f)\|_{\mu;p}\le c\|f\|_{\mu;p}$, and
\be\label{goodapprox}
E_{\mu;m,p}(f) \le \|f-\sigma_m(\nu,\mu; h,f)\|_{\mu;p} \le cE_{\mu;m/2,p}(f).
\ee
{\rm (c)} If $f$ is supported on a subinterval $I$ of $[-1,1]$, and $J$ is an interval with $I \subset J\subseteq [-1,1]$, then
\be\label{sigmalocalbd}
\|\sigma_m(\nu,\mu; h,f)\|_{\mu;\infty, [-1,1]\setminus J}\le c\|f\|_{\mu;1}m^{-S},
\ee
where $c$ may depend upon $I$ and $J$ in addition to $\mu$, $S$, and $h$.
\end{prop}

We point out an application of this proposition for the local information problem. Suppose $J$ is a subinterval of $[-1,1]$, $\C$ is a finite subset of $J$, and we have to construct an approximation to $f\in X^\infty(\mu)$ based on the values $\{f(z)\}_{z\in\C}$. We may take a subinterval $I$ of $J$, and a $C^\infty$ function $\phi$ which is equal to $1$ on $I$ and $0$ outside $J$. The above proposition  shows how to approximate $f$ on $I$ while keeping the growth of the approximating polynomial under control. We  find an integer $m$ such that there 
exists an M--Z quadrature measure $\nu_m$ of order $3m/2$ on a set containing $\C$ as its subset, and $\nu_m\preceq_\infty \mu$. The existence of such measures is proved in \cite[Proposition~5.1]{expframe} under the assumption that $\mu(\{x\})=0$ for each $x\in [-1,1]$.  Then \eref{goodapprox} shows that
$$
\|f-\sigma_{m}(\nu_m,\mu; h,f\phi)\|_{\mu;\infty, I}\le cE_{\mu;m/2,\infty}(f\phi),
$$
while \eref{sigmalocalbd} and \eref{uniformbd} limit the growth of $\|\sigma_{m}(\nu_m,\mu; h,f\phi)\|_{\mu;\infty, [-1,1]\setminus I}$.

\noindent
\textsc{Proof of Proposition~\ref{sigmaprop}.}
If $P\in\Pi_{m/2}$, $x\in\RR$, then $\Phi_m(\mu;h,x,\circ)P\in \Pi_{3m/2-1}$. Since $h(t)=1$ if $t\le 1/2$, it is easy to verify using \eref{genquad} that
$$
\int_{-1}^1 \Phi_m(\mu;h,x,t)P(t)d\nu(t) =\int_{-1}^1 \Phi_m(\mu;h,x,t)P(t)d\mu(t)=P(x).
$$
This proves part (a).
 
In view of \eref{genmzineq} and \eref{multmaskcond}, we have
$$
\sup_{x\in [-1,1]}\int_{-1}^1|\Phi_m(\mu;h,x,t)|d|\nu|(t)\le c \sup_{x\in [-1,1]}
\int_{-1}^1 |\Phi_m(\mu;h,x,t)|d\mu(t)\le c.
$$
The estimate \eref{uniformbd} is now a simple consequence of the Riesz--Thorin interpolation theorem and the fact that $\nu\preceq_p \mu$ (cf. \cite[Lemma~4.1]{locjacobi}).

If $P\in\Pi_{m/2}$ is arbitrary, then part (a) and \eref{uniformbd} imply that
$$
E_{\mu;m,p}(f) \le \|f-\sigma_m(\nu,\mu; h,f)\|_{\mu;p}=\|f-P-\sigma_m(\nu,\mu; h,f-P)\|_{\mu;p}\le c\|f-P\|_{\mu;p}.
$$
This proves \eref{goodapprox}.

The estimate \eref{sigmalocalbd} is easy to deduce using \eref{localizationcond}.
\qed

\begin{cor}\label{mzmeascor}
Let $m\ge 0$ be an integer, $\nu$ be an M--Z quadrature measure of order $3m-1$, and $\nu\preceq_\infty\mu$. Then 
\be\label{polynormequiv}
\|P\|_{\nu;p}\sim\|P\|_{\mu;p}, \qquad P\in\Pi_m, \ 1\le p\le \infty.
\ee
\end{cor}
\begin{Proof}
If $f\in L^1(\mu)$, then \eref{genmzineq} and \eref{uniformbd} (applied with $\mu$ in place of $\nu$) imply that 
$$
\|\sigma_{2m}^*(\mu;h,f)\|_{\nu;1}\le c\|\sigma_{2m}^*(\mu;h,f)\|_{\mu;1}\le c\|f\|_{\mu;1}.
$$
We note that $\nu\preceq_\infty \mu$ and that $\sigma_{2m}^*(\mu;h,f)\in C([-1,1])$. An application of \eref{uniformbd} again with $p=\infty$ and $\mu$ in place of $\nu$ shows that
$$
\|\sigma_{2m}^*(\mu;h,f)\|_{\nu;\infty}\le c\|\sigma_{2m}^*(\mu;h,f)\|_{\mu;\infty}\le c\|f\|_{\mu;\infty}, \qquad f\in L^\infty(\mu).
$$
Consequently, the Riesz--Thorin interpolation theorem yields $\|\sigma_{2m}^*(\mu;h,f)\|_{\nu;p}\le c\|f\|_{\mu;p}$ for every $p\in (1,\infty)$ and $f\in L^p(\mu)$. Using this estimate with $P\in\Pi_m$ in place of $f$, we obtain from part (a) of Proposition~\ref{sigmaprop} that $\|P\|_{\nu;p}\le c\|P\|_{\mu;p}$, $1\le p\le \infty$. The estimates in the reverse direction follow directly by using \eref{uniformbd} with $P$ in place of $f$ and recalling part (a) of Proposition~\ref{sigmaprop} again.
\end{Proof}

The following theorem describes the local adaptivity property of the summability operators.

\begin{theorem}\label{sigmatheo}
Let  $1\le p\le\infty$, $f\in X^p(\mu)$, $x_0\in [-1,1]$, $0<\rho\le\infty$, $\gamma>0$, $S>\max(1,\gamma)$. For each integer $n\ge 0$, let  $\nu_n$ be an M--Z quadrature measure of order $(3/2)(2^n)-1$, and $\nu_n\preceq_p \mu$. Then
 $f\in B_{\mu;p,\rho,\gamma}(x_0)$ if and only if there exists a nondegenerate interval $I$ centered at $x_0$ such that
\be\label{sigmalocapprox}
\|f-\sigma_{2^n}(\nu_n,\mu; h,f)\|_{\mu;p, I}\in\seqb_{\rho,\gamma}.
\ee
\end{theorem}

\begin{Proof} 
Let $f\in B_{\mu;p,\rho,\gamma}(x_0)$, and $J$ be an interval  centered at $x_0$ such that $E_{\mu;2^n,p}(f\phi)\in\seqb_{\rho,\gamma}$ for every $C^\infty$ function $\phi$ supported on $J$. We take $I\subset I_1\subset J$ to be  intervals centered at $x_0$, and a $C^\infty$ function $\psi$ supported on $J$ and equal to $1$ on $I_1$. Then \eref{sigmalocalbd} leads to
$$
\|\sigma_{2^n}(\nu_n,\mu; h,(1-\psi) f)\|_{\mu;p,I}\le c(I,J,f)2^{-nS}.
$$
 Since $\psi(t)=1$ for $t\in I$, we conclude that
\begin{eqnarray*}
\|f-\sigma_{2^n}(\nu_n,\mu; h, f)\|_{\mu;p,I}&\le& \|\psi f-\sigma_{2^n}(\nu_n,\mu; h, \psi f)\|_{\mu;p,I} + \|\sigma_{2^n}(\nu_n,\mu; h,(1-\psi) f)\|_{\mu;p,I}\\
&\le& c(I,J,f)\left\{E_{\mu;2^{n-1},p}(\psi f)+2^{-nS}\right\}.
\end{eqnarray*}
This proves \eref{sigmalocapprox}.

Conversely, let \eref{sigmalocapprox} hold, and $\phi$ be a $C^\infty$ function  supported on $I$. The direct theorem of approximation theory shows that there exists $R\in\Pi_{2^n}$ such that $\|\phi-R\|_\infty\le c(\phi)2^{-nS}$. Therefore, using \eref{uniformbd} and the fact that $\nu_n\preceq_p \mu$, we derive that
\begin{eqnarray*}
E_{\mu;2^{n+1},p}(f\phi) &\le&\|f\phi-R\sigma_{2^n}(\nu_n,\mu; h, f)\|_{\mu;p}\\
&\le&\|(f-\sigma_{2^n}(\nu_n,\mu; h, f))\phi\|_{\mu;p} +\|(\phi-R)\sigma_{2^n}(\nu_n,\mu; h, f)\|_{\mu;p}\\
&\le& c(\phi)\left\{\|f-\sigma_{2^n}(\nu_n,\mu; h, f)\|_{\mu;p,I} +2^{-nS}\|f\|_{\mu;p}\right\}.
\end{eqnarray*}
In view of \eref{sigmalocapprox}, this implies that $f\in B_{\mu;p,\rho,\gamma}(x_0)$.
\end{Proof}

The discrete Hardy inequality \cite[Lemma~3.4, p.~27]{devlorbk} implies that if $\{a_j\}_{j=0}^\infty\in \seqb_{\rho,\gamma}$ then $\{\sum_{j=k}^\infty a_j\}_{k=0}^\infty\in \seqb_{\rho,\gamma}$ as well. Therefore, we can deduce the following corollary using  \eref{goodapprox} and \eref{sigmalocapprox}.

\begin{cor}\label{framecor}
With the set up as in Theorem~\ref{sigmatheo}, we have
\be\label{seriesrep1}
f=\sum_{n=0}^\infty \left(\sigma_{2^n}(\nu_n,\mu; h, f)-\sigma_{2^{n-1}}(\nu_{n-1},\mu; h, f)\right)
\ee
with convergence in the sense of $X^p(\mu)$. Moreover, $f\in B_{\mu;p,\rho,\gamma}(x_0)$ if and only if there exists a nondegenerate interval $I$ centered at $x_0$ such that
$$
\|\sigma_{2^n}(\nu_n,\mu; h, f)-\sigma_{2^{n-1}}(\nu_{n-1},\mu; h, f)\|_{\mu;p, I}\in\seqb_{\rho,\gamma}.
$$
\end{cor}

The frame properties of the operators $f\mapsto\sigma_{2^n}(\nu_n,\mu; h, f)-\sigma_{2^{n-1}}(\nu_{n-1},\mu; h, f)$ were described in \cite{locjacobi}. 
We will now describe a tight frame, and demonstrate its use in characterization of local Besov spaces, thereby achieving all the objectives listed in the introduction.

Let $g(t):=\sqrt{h(t)-h(2t)}$. Then $g$ is supported on $[1/4,1]$, We assume that $g$ is also a multiplier mask of order $S$.  We define
$$
\tau_n^*(\mu;h,f):=\sigma_{2^n}^*(\mu;g,f).
$$
In the following theorem, we can choose $\nu_n$ to be a discretely supported positive measure. In the case of the Jacobi weights, it may be the measure that associates with each zero $x_{k,c(2^n)}(\mu)$, $k=1,\cdots,c(2^n)$ of $p_{c(2^n)}$  the mass $\lambda_{c(2^n)}(\mu;x_{k,c(2^n)}(\mu))$ for a suitable integer constant $c$ (\cite[Theorem~25, p.~168]{nevai}).

\begin{theorem}\label{tightframetheo}
Let $1\le p\le\infty$, $f\in X^p(\mu)$, $x_0\in [-1,1]$, $0<\rho\le\infty$, $\gamma>0$, $S>\max(1,\gamma)$. For each integer $n\ge 0$, let  $\nu_n$ be a sequence of measures with  each $\nu_n$ satisfying
\be\label{tightquad}
\int_{-1}^1 Pd\nu_n =\int_{-1}^1 Pd\mu, \qquad P\in \Pi_{2^{n+1}-1}.
\ee
{\rm (a)} We have, with convergence in the sense of $X^p$,
\be\label{seriesrep}
f = \sum_{n=0}^\infty \int_{-1}^1 \tau_n^*(\mu;h,f,t)\Phi_{2^n}(\mu;g,\circ,t)d\nu_n(t). 
\ee
{\rm (b)} In the case when $p=2$,
\be\label{parsevalframe}
\|f\|_{\mu;2}^2=\sum_{n=0}^\infty \|\tau_n^*(\mu;h,f)\|_{\nu_n;2}^2.
\ee
{\rm (c)} Let $\nu_n$ be an M--Z measure of order $3(2^n)-1$ and $\nu_n\preceq_p \mu$. Then $f\in B_{\mu;p,\rho,\gamma}(x_0)$ if and only if there exists a nondegenerate interval $I$ centered at $x_0$ such that $\left\{\|\tau_n^*(\mu;h,f)\|_{\mu;p,I}\right\}\in\seqb_{\rho,\gamma}$, or equivalently,  there exists a nondegenerate interval $I$ centered at $x_0$ such that $\left\{\|\tau_n^*(\mu;h,f)\|_{\nu_n;p,I}\right\}\in\seqb_{\rho,\gamma}$.
\end{theorem}

We remark that by choosing $\nu_n$ to be a discretely supported positive measure, the representation \eref{seriesrep} is actually a double sum representation, similar to classical wavelet expansions. The equation \eref{parsevalframe} then appears as the analogue of the classical Parseval identity for trignometric Fourier series. Part (c) shows that the absolute values of the coefficients $(\tau_n^*(\mu;h,f,t))_{t\in\mbox{ supp }(\nu_n)}$ in the representation \eref{seriesrep} can be used to characterize the local Besov spaces. It should be noted, however, that the functions $\Phi_{2^n}(\mu;g,\circ,t)$, $t\in\mbox{ supp }(\nu_n)$, are not necessarily linearly independent.

\noindent\textsc{Proof of Theorem~\ref{tightframetheo}.} In view of \eref{tightquad},  a straightforward calculation using the orthonormality of the polynomials $p_k(\mu)$ shows that for $x\in \RR$,
\bea\label{pf3eqn2}
\int_{-1}^1 \tau_n^*(\mu;h,f,t)\Phi_{2^n}(\mu;g,x,t)d\nu_n(t)&=&\int_{-1}^1\!\! \tau_n^*(\mu;h,f,t)\Phi_{2^n}(\mu;g,x,t)d\mu(t)\nonumber\\
&=&\sigma_{2^n}^*(\mu;h,f,x)-\sigma_{2^{n-1}}^*(\mu;h,f,x).
\eea
The equation \eref{seriesrep} follows from \eref{seriesrep1}, used with each $\nu_n=\mu$. Similarly,
$$
\|\tau_n^*(\mu;h,f)\|_{\nu_n;2}^2=\|\tau_n^*(\mu;h,f)\|_{\mu;2}^2,
$$
and \eref{parsevalframe} follows by a straightforward calculation (cf. \cite[Theorem~3]{fasttour}.) Next, we prove part (c), where stronger assumptions on $\nu_n$ are made. Since $g(t)=0$ if $t\le 1/4$, and $g$ is a multiplier mask of order $S$, we have for any $P\in\Pi_{2^{n-2}}$,
$$
\|\tau_n^*(\mu;h,f)\|_{\mu;p}=\|\tau_n^*(\mu;h,f-P)\|_{\mu;p}\le c\|f-P\|_{\mu;p}.
$$
Thus, $\|\tau_n^*(\mu;h,f)\|_{\mu;p}\le cE_{\mu;2^{n-2},p}(f)$. Therefore, if $f\in B_{\mu;p,\rho,\gamma}(x_0)$, then an argument similar to the proof of Theorem~\ref{sigmatheo}, using the assumption that $g$ is a multiplier mask of order $S$, implies that $\left\{\|\tau_n^*(\mu;h,f)\|_{\mu;p,I}\right\}\in\seqb_{\rho,\gamma}$ for some interval $I$ containing $x_0$. Conversely, let $\left\{\|\tau_n^*(\mu;h,f)\|_{\mu;p,I}\right\}\in\seqb_{\rho,\gamma}$, $J\subset I$ be a proper subinterval centered at $x_0$, and $\phi$ be an arbitrary $C^\infty$ function supported on $J$. Denoting by ${\mathcal X}$ the characteristic function of $I$, we obtain in view of \eref{multmaskcond} with $g$ in place of $h$ that 
\bea\label{pf3eqn1}
\lefteqn{\left\|\int_I \tau_n^*(\mu;h,f,t)\Phi_{2^n}(\mu;g,\circ,t)d\mu(t)\right\|_{\mu;p,J}}\nonumber\\
&\le& \left\|\int_{-1}^1 \tau_n^*(\mu;h,f,t){\mathcal X}(t)\Phi_{2^n}(\mu;g,\circ,t)d\mu(t)\right\|_{\mu;p}\le c\|\tau_n^*(\mu;h,f)\|_{\mu;p,I}.
\eea
If $x\in J$, then
$$
\int_{[-1,1]\setminus I}|\tau_n^*(\mu;h,f,t)\Phi_{2^n}(\mu;g,x,t)|d\mu(t)\le c(I,J)2^{-nS}\|\tau_n^*(\mu;h,f)\|_{\mu;p}\le c(I,J)2^{-nS}\|f\|_{\mu;p}.
$$
The above two estimates and \eref{pf3eqn2} show that
\begin{eqnarray*}
\|\sigma_{2^n}^*(\mu;h,f)-\sigma_{2^{n-1}}^*(\mu;h,f)\|_{\mu;p,J} &=&\left\|\int_{-1}^1 \tau_n^*(\mu;h,f,t)\Phi_{2^n}(\mu;g,\circ,t)d\mu(t)\right\|_{\mu;p,J}\\
&\le&c_1\|\tau_n^*(\mu;h,f)\|_{\mu;p,I} +c(I,J)2^{-nS}\|f\|_{\mu;p}.
\end{eqnarray*}
Thus, $\|\sigma_{2^n}^*(\mu;h,f)-\sigma_{2^{n-1}}^*(\mu;h,f)\|_{\mu;p,J}\in\seqb_{\rho,\gamma}$, and Corollary~\ref{framecor} implies that $f\in B_{\mu;p,\rho,\gamma}(x_0)$. \qed

\bhag{Quadrature formulas}\label{quadsect}
 In this section, we restrict ourselves  to the case of Jacobi weights, and prove the existence of quadrature formulas at arbitrarily chosen points on the interval, exact for integrating high degree polynomials with respect to the Jacobi weights.

For $\a,\beta>-1$, let
$w_{\a,\beta}(x):=(1-x)^\a(1+x)^\beta$, $x\in (-1,1)$, and $w_{\a,\beta}(x)=0$ otherwise. We define $d\mu^{(\alpha,\beta)}(x)=w_{\a,\beta}(x)dx$, and denote the corresponding degree $k$ orthonormalized polynomial $p_k(\mu^{(\alpha,\beta)})$ with positive leading coefficient by $p_k^{(\alpha,\beta)}$, ($k=0,1,\cdots$). In particular, for integers $k, j=0,1,\cdots$,
$$
\int_{-1}^1 \jsup{p_k}{\a}{\b}(x) \jsup{p_j}{\a}{\b}(x)w_{\a,\b}(x)dx=\left\{\begin{array}{ll}
1, &\mbox{ if $k=j$,}\\
0, &\mbox{otherwise.}
\end{array}\right.
$$
In the case $\alpha=\beta=-1/2$, one obtains the Chebyshev polynomials. For $\theta\in [0,\pi]$, let $T_k(\cos\theta):=\cos k\theta$. We have for $x\in [-1,1]$,
$p_0^{(-1/2,-1/2)}(x)=(1/\sqrt\pi)T_0(x)$, and for $k=1,2,\cdots,$ $p_k^{(-1/2,-1/2)}(x)=(\sqrt{2/\pi})T_k(x)$. 

Let $M\ge 4$ be an integer, $0=\theta_0\le \theta_1<\cdots<\theta_k<\theta_{k+1}<\cdots<\theta_M\le \theta_{M+1}=\pi$ be arbitrary points, and $z_k=\cos\theta_k$, $k=0,\cdots,M+1$. We will tacitly assume that the set $\C=\{z_k\}$ is one of the members of a nested sequence of sets of points on $[-1,1]$, whose union in dense in $[-1,1]$. Thus, all constants will be independent of $M$ and the points $z_k$. The \emph{mesh norm} (also known as fill distance) of the set is defined by
$$
\delta_{\C}:= \max_{\theta\in [0,\pi]}\min_{1\le k\le M} |\theta-\theta_k|=\frac{1}{2}\max_{0\le k\le M} |\theta_{k+1}-\theta_k|.
$$
The \emph{separation radius} of $\C$ is defined by $q_\C:=(1/2)\min_{0\le k\le M} |\theta_{k+1}-\theta_k|$. The set $\C$ is called $\rho$-uniform if $\delta_\C\le \rho q_\C$, and uniform if it is $\rho$--uniform for some $\rho$. Again, it is to be understood that the value of $\rho$ is the same for the whole implicitly understood sequence of sets $\C$.  If $\delta_\C\le \pi/4$, and $m\ge 1$ is the integer part of $\pi\delta_\C^{-1}/4$, then each interval of the form $\disp [\frac{(4j+1)\pi}{4m}, \frac{(4j+3)\pi}{4m}]$ contains at least one $\theta_k$. We choose one $\theta_{k_j}$ from each such interval to obtain a subset $\C'=\{\cos\theta_{k_j}\}\subset\C$ such that $q_{\C'}\ge \pi/(4m)$ and $\delta_\C\le \delta_{\C'} \le \pi/m \le 4\delta_\C$.  Moreover, $\{\theta_{k_j}\}\subset [c/m, \pi-c/m]$. Thus, in the sequel, we may assume that the sets $\C$ are all $4$-uniform, and that there exists a constant $c$ with the property that $\{\cos^{-1}z\ :\ z\in\C\}\subset [c/m, \pi-c/m]$. Clearly, for a uniform set $\C$, $q_\C\sim\delta_\C\sim |\C|^{-1}$.

\begin{theorem}\label{jacobiquadtheo}
Let $\a,\beta>-1$, $\rho\ge 1$. There exists a constant $a>0$ depending only on $\a,\beta,\rho$, with the following property. Let $\C\subset (-1,1)$ be a finite $\rho$--uniform set, and $m\ge 1$ be an integer such that $\pi/(2m)\le \delta_\C\le \pi/m$. Then there exist positive numbers $w_z$, $z\in\C$, with the following properties:
\be\label{quadformula}
\sum_{z\in \C} w_zP(z) =\int_{-1}^1 P(t)w_{\alpha,\beta}(t)dt, \qquad P\in \Pi_{am}, 
\ee
\be\label{wtbds}
w_z\sim \lambda_m(\jsup{\mu}{\a}{\b};z), \qquad z\in\C,
\ee
and
\be\label{mzquadineq}
\sum_{z\in\C}w_z|P(z)| \sim \int |P(t)|w_{\alpha,\beta}(t)dt, \qquad P\in \Pi_{am}. 
\ee
\end{theorem}

The proof of this theorem follows a familiar theme, introduced in \cite{JetSW2, mnw1}. We first recall an abstract quadrature formula in the setting of general finite dimensional spaces (cf. \cite[Theorem~3.2.1]{indiapap}). 

\begin{prop}\label{abstractquadprop}
Let $(X,\|\cdot\|_X)$ be a finite dimensional normed linear space,  $(X^*,\|\cdot\|_{X^*})$ be its dual space, $\Z=\{x_1^*,\cdots,x_M^*\}\subset X^*\setminus\{0\}$,  and $x^*\in X^*$. Suppose the operator $x\mapsto 
 (x_1^*(x),\cdots,x_M^*(x))$, $x\in X$, is one-to-one, and the following two conditions are satisfied: (1) If $x\in X$ and $x_\ell^*(x) \ge 0$ for $\ell=1,\cdots, M$, then $x^*(x)\ge 0$, and (2) there exists some $x_0\in X$ such 
that $x^*_\ell(x_0) >0$ for $\ell=1,\cdots,M$. Then there exist nonnegative numbers $W_\ell$, $\ell=1,\cdots,M$ such that
\be\label{absquad}
x^*(x)=\sum_{\ell=1}^M W_\ell x_\ell^*(x).
\ee
\end{prop}

We will apply this proposition with $X=\Pi_{am}$ for a judiciously chosen $a$, and $x_\ell^*(P)=P(x_\ell)$ for $x_\ell\in\C$. The choice of $a$ is dictated by the following Lemma~\ref{mzlemma}.
The inequalities of the form \eref{mzineq} are known as Marcinkiewicz--Zygmund (M--Z) inequalities. 

\begin{lemma}\label{mzlemma}
Let $\a,\beta>-1$, $\rho\ge 1$. There exists a constant $a>0$ depending only on $\a,\beta,\rho$, with the following property. Let $\C\subset (-1,1)$ be a finite $\rho$--uniform set, and $m\ge 1$ be an integer such that $\pi/(2m)\le \delta_\C\le \pi/m$. Then for every $P\in \Pi_{am}$,
\be\label{mzineq}
c_1\|P\|_{\jsup{\mu}{\a}{\beta};1} \le \sum_{z\in\C} \lambda_m(\jsup{\mu}{\a}{\b};z)|P(z)|\le c_2\|P\|_{\jsup{\mu}{\a}{\beta};1}.
\ee
Moreover, if $P(z)\ge 0$ for each $z\in\C$, then there exists $c\in (0,1)$ (depending only on $\a,\b$)  such that
\be\label{positivity}
\int_{-1}^1 P(t)w_{\alpha,\beta}(t)dt\ge c\sum_{z\in\C} \lambda_m(\jsup{\mu}{\a}{\b};z)P(z)\ge 0.
\ee
\end{lemma}

In order to prove this lemma, we need to recall certain basic inequalities. Parts (a) and (b) can be found in \cite[Theorem~4]{mastr} and part (c) can be found in  \cite[p.~108]{nevai}. We recall that we omit the mention of the measure $\mu$ from the notations if it is the Lebesgue measure.

\begin{prop}\label{ineqprop}
Let $m\ge 1$ be an integer, $P\in\Pi_m$, $1\le p\le\infty$, $\a,\beta> -1$, and $w_{\a,\beta}\in L^p$. Let 
\be\label{barjacobiwt}
\overline{w}_{m,\a,\beta}(x):=(\sqrt{1-x}+1/m)^{2\a}(\sqrt{1+x}+1/m)^{2\beta}.
\ee
{\rm (a)} {\rm (Markov--Bernstein inequality)} 
\be\label{markovineq}
\|P'\overline{w}_{m,\a+1/2,\beta+1/2}\|_p \le cm\|Pw_{\a,\beta}\|_p.
\ee
{\rm (b)} {\rm (Nikolskii inequality)}
\be\label{nikolskii}
\|P\overline{w}_{m,\a,\beta}\|_r\le cm^{2(1/r-1/p)}\|Pw_{\a,\beta}\|_p, \qquad 1\le p\le r\le\infty.
\ee
{\rm (c)}
 Let $\theta,\varphi\in [0,\pi]$, $m\ge 1$ be an integer, $|\theta-\varphi|\sim 1/m$, $x=\cos\theta$, $y=\cos\varphi$, $x\le y$ and $u,v\in [x,y]$. Then 
\be\label{basicfact}
\lambda_m(\jsup{\mu}{\a}{\b};u)\sim (1/m)\overline{w}_{m,\alpha+1/2,\beta+1/2}(u)\sim (1/m)\overline{w}_{m,\alpha+1/2,\beta+1/2}(v).
\ee
\end{prop}

\noindent\textsc{Proof of Lemma~\ref{mzlemma}.}
In this proof only, let  $M\ge 4$ be an integer, $\C=\{z_k=\cos\theta_k\}_{k=1}^M$. Since $\C$ is uniform, we also have $q_\C\sim 1/m$. In this proof  only, let $\tilde\theta_0=0$, $\tilde\theta_j=\theta_j+q_\C$, $j=1,\cdots,M-1$, $\tilde\theta_M=\pi$, $\tilde z_j=\cos\tilde\theta_j$, $0\le j\le M$, and $I_j:=[\tilde z_j,\tilde z_{j-1}]$, $j=1,\cdots,M$. Let $n\ge 1$ be an integer, and $P\in\Pi_n$ be arbitrary. We will prove first that
\be\label{basicmzineq}
\sum_{j=1}^{M}\int_{I_j}|P(t)-P(z_{j})|w_{\a,\beta}(t)dt \le c_3\left\{\frac{n}{m}+\left(\frac{n}{m}\right)^2\right\}\|Pw_{\a,\beta}\|_1. 
\ee

We observe that $1-\tilde z_1=2\sin^2(\tilde\theta_1/2) \sim 1/m^2$. In view of \eref{basicfact},
\be\label{pf1eqn3}
\int_{\tilde z_1}^1\!\! w_{\a,\beta}(t)dt\sim (1-\tilde z_1)^{\a+1}\sim \frac{\overline{w}_{m,\a,\beta}(\tilde z_1)}{m^2}\sim \frac{\overline{w}_{m,\a,\beta}(z_1)}{m^2}\sim \frac{\overline{w}_{m,\a+1/2,\beta+1/2}(z_1)}{m}\sim\lambda_m(\jsup{\mu}{\a}{\b};z_1).
\ee
Similarly, 
\bea\label{pf1eqn4}
\int_{-1}^{\tilde z_{M-1}}\!\! w_{\a,\beta}(t)dt &\sim& (1+\tilde z_{M-1})^{\beta+1}\sim \frac{\overline{w}_{m,\a,\beta}(\tilde z_{M-1})}{m^2}\sim \frac{\overline{w}_{m,\a,\beta}(z_M)}{m^2}\nonumber\\
&\sim& \frac{\overline{w}_{m,\a+1/2,\beta+1/2}(z_M)}{m}\sim\lambda_m(\jsup{\mu}{\a}{\b};z_M).
\eea
In view of \eref{pf1eqn3} and \eref{nikolskii}, we have for $u\in [\tilde z_1,1]$,
$$
|P(u)|\int_{\tilde z_1}^1\!\! w_{\a,\beta}(t)dt\le   \frac{c}{m^2}|P(u)|\overline{w}_{m,\a,\beta}(\tilde z_1)  \le \frac{c}{m^2}|P(u)|\overline{w}_{m,\a,\beta}(u) \le \frac{cn^2}{m^2}\|Pw_{\a,\beta}\|_1.
$$
We estimate $|P(v)|\int_{-1}^{\tilde z_{M-1}}\!\! w_{\a,\beta}(t)dt$ ($v\in [-1,\tilde z_{M-1}]$) in the same way using \eref{pf1eqn4} in place of \eref{pf1eqn3}, and deduce that
\be\label{pf1eqn1}
\int_{\tilde z_1}^1\!\! |P(t)-P(z_1)|w_{\a,\beta}(t)dt+ \int_{-1}^{\tilde z_{M-1}}\!\!|P(t)-P(z_M)| w_{\a,\beta}(t)dt\le  \frac{cn^2}{m^2}\|Pw_{\a,\beta}\|_1. 
\ee
If $2\le j\le M-1$, then for $t, u\in I_j$, we have 
$$
w_{\a+1/2,\beta+1/2}(t)\sim \overline{w}_{m,\a+1/2,\beta+1/2}(t) \sim w_{\a+1/2,\beta+1/2}(u).
$$
 Consequently, using \eref{markovineq}, we obtain
\bea\label{pf1eqn2}
\lefteqn{\sum_{j=1}^{M-1}\int_{I_j}|P(t)-P(z_{j})|w_{\a,\beta}(t)dt  \le   \sum_{j=2}^{M-2}\int_{I_j}\int_{I_j}|P'(u)|du w_{\a,\beta}(t)dt}\nonumber\\
&\le& c \sum_{j=1}^{M-1}\int_{I_j}\int_{I_j}|P'(u)|w_{\a+1/2,\beta+1/2}(u)du (1-t^2)^{-1/2}dt \nonumber\\
&\le& c\max (\tilde \theta_{j+1}-\tilde\theta_j)\int_{-1}^1|P'(u)|\overline{w}_{m,\a+1/2,\beta+1/2}(u)du\nonumber\\
&\le& \frac{cn}{m}\|Pw_{\a,\beta}\|_1.
\eea
 The estimates \eref{pf1eqn1} and \eref{pf1eqn2} together imply \eref{basicmzineq}.

In the remainder of this proof, we consider the value of $c_3$ fixed as in \eref{basicmzineq}.  Let $a=\min (1, 1/(8c_3))$, and $P\in\Pi_{am}$. Then with $am$ in place of $n$,
\be\label{pf1eqn7}
c_3\left\{\frac{n}{m}+\left(\frac{n}{m}\right)^2\right\}\le c_3(a+a^2)\le 2c_3a\le  1/4.
\ee
Since 
$$
\left|\int_{-1}^1 |P(t)|w_{\a,\beta}(t)dt -\sum_{j=1}^M |P(z_j)|\int_{I_j}w_{\a,\beta}(t)dt\right|\le \sum_{j=1}^{M}\int_{I_j}|P(t)-P(z_{j})|w_{\a,\beta}(t)dt,
$$
we obtain from \eref{basicmzineq} that
\be\label{pf1eqn5}
(3/4)\|Pw_{\a,\beta}\|_1 \le \sum_{j=1}^M |P(z_j)|\int_{I_j}w_{\a,\beta}(t)dt\le (5/4)\|Pw_{\a,\beta}\|_1.
\ee
We now observe that for $2\le j\le M-1$, and $t\in I_j$, $w_{\a,\beta}(t)\sqrt{1-t^2}\sim \overline{w}_{m,\a+1/2,\beta+1/2}(t)$. Therefore, \eref{basicfact} implies that 
$$
\int_{I_j}w_{\a,\beta}(t)dt \sim \int_{I_j}\overline{w}_{m,\a+1/2,\beta+1/2}(t)\frac{dt}{\sqrt{1-t^2}}\sim \frac{\overline{w}_{m,\a+1/2,\beta+1/2}(z_j)}{m}\sim \lambda_m(\jsup{\mu}{\a}{\b};z_j).
$$
 Together with \eref{pf1eqn3} and \eref{pf1eqn4}, we have proved that
\be\label{pf1eqn6}
\int_{I_j}w_{\a,\beta}(t)dt \sim \lambda_m(\jsup{\mu}{\a}{\b};z_j), \qquad j=1,\cdots,M.
\ee
The estimates \eref{pf1eqn5} and \eref{pf1eqn6} imply \eref{mzineq}.

Next, let $P(z_j)\ge 0$ for $j=1,\cdots,M$. In view of \eref{basicmzineq} and \eref{pf1eqn5}, we obtain
\begin{eqnarray*}
\lefteqn{\left|\int_{-1}^1 P(t)w_{\a,\beta}(t)dt -\sum_{j=1}^M P(z_j)\int_{I_j}w_{\a,\beta}(t)dt\right|}\\
&\le& \sum_{j=1}^{M}\int_{I_j}|P(t)-P(z_{j})|w_{\a,\beta}(t)dt\le \frac{1}{4}\|Pw_{\a,\beta}\|_1\le \frac{1}{3}\sum_{j=1}^M P(z_j)\int_{I_j}w_{\a,\beta}(t)dt.
\end{eqnarray*}
Consequently,
$$
\int_{-1}^1 P(t)w_{\a,\beta}(t)dt\ge \frac{2}{3}\sum_{j=1}^M P(z_j)\int_{I_j}w_{\a,\beta}(t)dt\ge 0.
$$
In view of \eref{pf1eqn6}, this implies \eref{positivity}.
\qed

We are now in a position to prove Theorem~\ref{jacobiquadtheo}, mainly using the ideas in \cite{mnw1}, but using a trick in \cite{narcpetward} to prove \eref{wtbds}.

\noindent
\textsc{Proof of Theorem~\ref{jacobiquadtheo}.} We let $a$ be as in Lemma~\ref{mzlemma}. In Proposition~\ref{abstractquadprop}, we  choose $\Pi_{am}$ in place of $X$, $|\C|$ in place of $M$, the mappings $P\mapsto P(z)$, $z\in \C$, as the set $\Z$ of linear functionals. The estimate \eref{mzineq} shows that the operator $P\mapsto (P(z))_{z\in\C}$ is one-to-one. We will take 
$$
x^*(P)=\int_{-1}^1 P(t)d\jsup{\mu}{\a}{\beta}(t) -\frac{c}{2}\sum_{z\in\C} \lambda_m(\jsup{\mu}{\a}{\b};z)P(z),
$$
where $c$ is the constant appearing in \eref{positivity}. 
If each $P(z)\ge 0$, $z\in\C$, then \eref{positivity} implies that $x^*(P)\ge 0$. Taking $x_0$ in Proposition~\ref{abstractquadprop} to be the polynomial identically equal to $1$, we see from \eref{positivity} that $x^*(x_0)>0$. Thus, all the conditions of Proposition~\ref{abstractquadprop} are satisfied, and we obtain nonnegative $W_z$, $z\in\C$, such that
$$
\int_{-1}^1 P(t)w_{\a,\beta}(t)= \sum_{z\in \C}(W_z+c\lambda_m(\jsup{\mu}{\a}{\b};z)/2)P(z), \qquad P\in \Pi_{am}.
$$
Setting $w_z= W_z+c\lambda_m(\jsup{\mu}{\a}{\b};z)/2$, we have proved \eref{quadformula}, and also that $w_z\ge c\lambda_m(\jsup{\mu}{\a}{\b};z)/2$, $z\in\C$.

Next, let $\xi\in \C$, and $n$ be the integer part of $am/2$. Then \eref{quadformula} and \eref{darbouxreprod} imply that
$$
w_\xi K_n^2(\jsup{\mu}{\a}{\beta}; \xi,\xi) \le \sum_{z\in\C}w_zK_n^2(\jsup{\mu}{\a}{\beta}; \xi,z)=\int_{-1}^1 K_n^2(\jsup{\mu}{\a}{\beta}; \xi,t)w_{\a,\beta}(t)=K_n(\jsup{\mu}{\a}{\beta}; \xi,\xi).
$$
Hence, 
\be\label{pf2eqn2}
w_\xi \le \lambda_n(\jsup{\mu}{\a}{\b};\xi), \qquad \xi\in\C.
\ee
In view of \eref{basicfact} and the already proved fact that $w_\xi\ge (c/2)\lambda_m(\jsup{\mu}{\a}{\b};\xi)$, $\xi\in\C$, this implies \eref{wtbds}. 

The estimates \eref{mzquadineq} are clear from \eref{wtbds} and \eref{mzineq}.
\qed

\noindent
\textbf{Remark.} We have actually proved \eref{pf2eqn2} for any positive numbers $w_z$ for which \eref{quadformula} is valid for all $P\in\Pi_n$, whether obtained via Theorem~\ref{jacobiquadtheo} or not.

\bhag{Localized polynomial kernels}\label{kernsect}
The objective of this section is to demonstrate that in the case of the Jacobi weights, any sufficiently smooth, compactly supported function $h$, which is constant on a neighborhood of $0$, is a summability mask of order $S$. This fact was essentially proved in \cite{locjacobi}, but our estimates here are more elegant. We prefer to generalize the summability kernels defined in Section~\ref{framesect} to define
\be\label{phidef}
\Phi(\mu; {\bf h},x,t):=\sum_{k=0}^\infty h_kp_k(\mu;x)p_k(\mu;t), \qquad x,t\in\RR,
\ee
where ${\bf h}=\{h_k\}_{k=0}^\infty$ is a compactly supported sequence. 

The  summability operator corresponding to the kernel $\Phi(\mu; {\bf h})$ is defined by
\be\label{sigmastardef}
\sigma^*(\mu;{\bf h},f,x):=\int_{-1}^1 \Phi(\mu;{\bf h},x,t)f(t)d\mu(t)=\sum_{k=0}^\infty h_k\hat f(\mu;k)p_k(\mu;x).
\ee
  An interesting problem in harmonic analysis is to find conditions on the sequence ${\bf h}$ so that this operator is bounded in some $L^p(\mu)$; i.e.,
$\|\sigma^*(\mu;{\bf h}, f)\|_{\mu;p}\le c({\bf h})\|f\|_{\mu;p}$. If this is the case, the sequence ${\bf h}$ is called a multiplier sequence. In the important case when $\mu=\jsup{\mu}{\a}{\beta}$, several conditions to ensure that ${\bf h}$ is a multiplier sequence  are available in the literature (e.g. \cite{muckenhoupt, gasper}). Typically, these conditions do not hold when $p=1$ or $p=\infty$. For our research, we therefore need to find some conditions to ensure that the sequence ${\bf h}$ is a multiplier sequence also in these cases. 
It is not difficult to see (cf. \cite[Proposition~2.1]{mauropap}) that this requirement is equivalent to 
\be\label{generalmultiplier}
\sup_{x\in [-1,1]}\int_{-1}^1 |\Phi(\mu; {\bf h},x,t)|d\mu(t) \le c({\bf h}).
\ee

\begin{theorem}\label{jacobikerntheo}
Let $\alpha,\beta\ge -1/2$, $S\ge 1$ be an integer, $h_k=0$ for all sufficiently large $k$. Then for $\theta,\varphi\in [0,\pi]$,
\bea\label{jacobikernestgen}
\lefteqn{|\Phi(\jsup{\mu}{\a}{\b};{\bf h},\cos\theta,\cos\varphi)|}\nonumber\\
&\le& c\left\{\begin{array}{ll}
\disp\sum_{k=0}^\infty \min\left((k+1), \frac{1}{|\theta-\varphi|}\right)^{\max(\alpha,\beta)+S+1/2}\times\\
\disp\qquad\times\sum_{m=1}^{S}(k+1)^{\max(\alpha,\beta)+1/2-S+m}|\Delta^{m}h_k|.
\end{array}\right.
\eea
In particular, if $h :[0,\infty)\to [0,\infty)$ is a compactly supported function that can be expressed as an $S$ times iterated integral of a function of bounded variation, and $h'(t)=0$ in a neighborhood of $0$, then for $\lambda \ge 1$,
\be\label{jacobikernestmain}
|\Phi_\lambda(\jsup{\mu}{\a}{\b};h,\cos\theta,\cos\varphi)|\le c\lambda^{2\max(\a,\beta)+2}V(\derf{h}{S-1})\min\left(1, \frac{1}{(\lambda|\theta-\varphi|)^{\max(\alpha,\beta)+S+1/2}}\right).
\ee
\end{theorem}

The following corollary is a routine consequence of Theorem~\ref{jacobikerntheo}, proved in \cite[Lemma~4.6]{locjacobi} using a different argument. We omit the proof.

\begin{cor}\label{jacobicor}
Let $\alpha,\beta\ge -1/2$, $S\ge \max(\a,\beta)+3/2$ be an integer, $h_k=0$ for all sufficiently large $k$. Then
\be\label{multpliercond}
\sup_{ x\in [-1,1]}\int_{-1}^1 |\Phi(\jsup{\mu}{\a}{\b};{\bf h},x,y)|d\jsup{\mu}{\alpha}{\beta}(y) \le c\sum_{j=1}^{S}\sum_{k=0}^\infty  (k+1)^{j-1}|\Delta^j h_k|.
\ee
In particular, if $h :[0,\infty)\to [0,\infty)$ is a compactly supported function that can be expressed as an $S$ times iterated integral of a function of bounded variation,  then
\be\label{jacobimultmask}
\sup_{\lambda>0,\ x\in [-1,1]}\int_{-1}^1 |\Phi_\lambda(\jsup{\mu}{\a}{\b};h,x,y)|d\jsup{\mu}{\alpha}{\beta}(y) <\infty.
\ee
\end{cor}
Analogues of \eref{jacobikernestmain} and \eref{jacobimultmask} hold in very general situations, e.g., in the case of eigenfunctions of Laplace--Beltrami operators on a smooth manifold rather than Jacobi polynomials \cite{mauropap}.

Although \eref{jacobimultmask} is routinely proved using an estimate similar to \eref{jacobikernestmain}, it is sometimes possible to obtain such estimates without first proving a localization estimate (\cite{freudbook, freud2}). We have recently observed in \cite{expframe} that there is a simple construction of exponentially localized operators from kernels satisfying analogues of \eref{jacobimultmask}, with interesting consequences for spectral approximation of piecewise analytic functions.

In order to prove Theorem~\ref{jacobikerntheo}, we start by recalling the following estimate proved in
 \cite[Lemma~4.10]{locjacobi}, using the explicit formulas for the Christoffel--Darboux kernel for the Jacobi polynomials: 

\begin{prop}\label{jacobikernprop}
Let $\alpha,\beta\ge -1/2$, $S\ge 1$ be an integer, $h_k=0$ for all sufficiently large $k$. Then
\bea\label{endptests1}
\lefteqn{|\Phi(\jsup{\mu}{\a}{\b};{\bf h},1,y)|}\nonumber\\
&\le& c\cases{\disp\sum_{k=0}^\infty \min\left((k+1)^2, \frac{1}{1-y}\right)^{\alpha/2+S/2+1/4}\times\cr
\disp\qquad\times\sum_{m=0}^{S-1}(k+1)^{\alpha+1/2-m}|\Delta^{S-m}h_k|, 
& if $0\le y<1$,\cr
\disp\sum_{k=0}^\infty (k+1)^{\alpha+\beta+1} \sum_{m=0}^{S-1}(k+1)^{-m}|\Delta^{S-m}h_k|, 
& if $-1\le y <0$,\cr}\nonumber\\
&&
\eea
and 
\bea\label{endptests2}
\lefteqn{|\Phi(\jsup{\mu}{\a}{\b};{\bf h},-1,y)|}\nonumber\\
&\le& c\cases{\disp\sum_{k=0}^\infty \min\left((k+1)^2, \frac{1}{1+y}\right)^{\beta/2+S/2+1/4}\times\cr
\qquad \disp\times\sum_{m=0}^{S-1}(k+1)^{\beta+1/2-m}|\Delta^{S-m}h_k|, 
& if $-1< y \le 0$,\cr
\disp\sum_{k=0}^\infty (k+1)^{\alpha+\beta+1} \sum_{m=0}^{S-1}(k+1)^{-m}|\Delta^{S-m}h_k|, 
& if $0< y\le1$.\cr}\nonumber\\
&&
\eea
\end{prop}
 
\noindent\textsc{Proof of Theorem~\ref{jacobikerntheo}.} We note first that the uniqueness of orthonormalized polynomials implies that $\jsup{p_k}{\a}{\beta}(x)=(-1)^k\jsup{p_k}{\beta}{\alpha}(-x)$, and hence, 
$\Phi(\jsup{\mu}{\alpha}{\beta}; {\bf h}, x,y)= \Phi(\jsup{\mu}{\beta}{\alpha}; {\bf h}, -x,-y)$. Therefore, we may assume that $\a\ge \beta$. The estimate \eref{endptests1} is then equivalent to \eref{jacobikernestgen} for the case when $\cos\theta=1$. 
To extend the estimate for $\Phi(\jsup{\mu}{\a}{\b};{\bf h},x,y)$ for every $x,y\in [-1,1]$,   we recall the following product formula \eref{prodformula} by Koornwinder \cite{koorn}: 
Let $\alpha\ge\beta\ge-1/2$, ${\cal R}:=[0,1]\times [0,\pi]$, and for $x,y\in [-1,1]$, $r\in [0,1]$, $\omega\in [0,\pi]$, let 
\be\label{koornfdef}
F(x,y;r,\omega):= \frac{(1+x)(1+y)}{2} + \frac{(1-x)(1-y)}{2}r^2 + \sqrt{1-x^2}\sqrt{1-y^2}r\cos\omega -1.
\ee
There exists a probability measure $\tilde\mu=\tilde\mu_{\alpha,\beta}$ on ${\cal R}$ such that  for $n=0,1,\cdots$, and $x,y\in [-1,1]$,
\be\label{prodformula}
\jsup{p_n}{\alpha}{\beta}(x)\jsup{p_n}{\alpha}{\beta}(y)=\int_{{\cal R}}\jsup{p_n}{\alpha}{\beta}(1)\jsup{p_n}{\alpha}{\beta}(F(x,y;r,\omega))d\tilde\mu(r,\omega).
\ee 
It follows that
\be\label{pf4eqn1}
\Phi(\jsup{\mu}{\alpha}{\beta}; {\bf h}, x,y)=\int_{{\cal R}}\Phi(\jsup{\mu}{\alpha}{\beta}; {\bf h}, 1,F(x,y;r,\omega))d\tilde\mu(r,\omega).
\ee
 Let $x= \cos \theta, y= \cos \varphi$. Then
\begin{eqnarray*}
\lefteqn{F(x,y;r ,\omega) = \frac{1}{2} (1+x)(1+y) + \frac{1}{2} (1-x)(1-y) + \sqrt{1-x^2} \sqrt{1-y^2}}\\
&&\qquad + \frac{1}{2} (1-x)(1-y)(r^2-1) +\sqrt{1-x^2} \sqrt{1-y^2} (r \cos \omega -1) -1\\
&=& xy+\sqrt{1-x^2} \sqrt{1-y^2} + \frac{1}{2}(1-x)(1-y)(r^2-1) + \sqrt{1-x^2} \sqrt{1-y^2} (r \cos \omega -1),
\end{eqnarray*}
and we have
\begin{eqnarray*}
 \lefteqn{1-F(x,y;r, \omega)=
 1-(xy+ \sqrt{1-x^2} \sqrt{1-y^2}) + \frac{1}{2} (1-x)(1-y)(1-r^2)}\\
&&\qquad + (1-r \cos \omega) \sqrt{1-x^2} \sqrt{1-y^2}\\
&\geq& 1- (xy+ \sqrt{1-x^2} \sqrt{1-y^2}) = 1 - \cos(\theta-\varphi) = 2 \sin^2 \frac{\theta- \varphi}{2}\\
& \geq& \frac{2}{\pi^2} (\theta-\varphi)^2.
\end{eqnarray*}
We now observe that $\a+S+1/2\ge \beta+1/2$, so that  
\eref{endptests1} implies
\begin{eqnarray*}
\lefteqn{|\Phi(\jsup{\mu}{\a}{\b};{\bf h},1,F(x,y;r, \omega))|}\\
&\le&c\disp\sum_{k=0}^\infty \min\left((k+1)^2, \frac{1}{1-F(x,y;r, \omega)}\right)^{\alpha/2+S/2+1/4}\times\\
&&\qquad\qquad\disp\qquad\times\sum_{m=0}^{S-1}(k+1)^{\alpha+1/2-m}|\Delta^{S-m}h_k|\\
&\le&c\disp\sum_{k=0}^\infty \min\left((k+1)^2, \frac{1}{(\theta-\varphi)^2}\right)^{\alpha/2+S/2+1/4}\times\\
&&\qquad\qquad\disp\qquad\times\sum_{m=0}^{S-1}(k+1)^{\alpha+1/2-m}|\Delta^{S-m}h_k|.
\end{eqnarray*}
Since $\tilde\mu$ is a probability measure, \eref{pf4eqn1} now implies \eref{jacobikernestgen}.

We observe next that  if $h :[0,\infty)\to [0,\infty)$ is a compactly supported function that can be expressed as an $S$ times iterated integral of a function of bounded variation,  $h'(t)=0$ in a neighborhood of $0$, and  $\lambda \ge 1$, then we may apply the above estimates with the sequence ${\bf h}_\lambda$ whose $k$-th term is given by $h_{\lambda,k}=h(k/\lambda)$, $k=0,1,\cdots$. A repeated application of mean value theorem implies that for any $s\in\RR$ and integer $r\ge 1$,
\be\label{hfuncdiffest}
\sum_{k=0}^\infty (k+1)^s\Delta^r h_{\lambda,k}\le c\lambda^{-s}\sum_{c\lambda\le k\le c_1\lambda} |\Delta^r h_{k,\lambda}|\le c\lambda^{s-r+1}V(h^{(r-1)}).
\ee
The estimate \eref{jacobikernestmain} follows from \eref{jacobikernestgen} and \eref{hfuncdiffest}. 
\qed


\bhag{Some comments on computations}\label{computesect}
There are several  algorithms for computing orthogonal polynomials and the corresponding Gauss--Jacobi quadrature formulas, as well as expressions of the form $\sum_{k=0}^n a_kp_k(\mu;x)$ \cite{gautschi}. Although we have not necessarily used these, the localization properties of our operators and their ability to detect singularities and local Lipschitz exponents has been demonstrated in a number of papers, in particular, \cite{loctrigwav, expframe, quadconst}.  We offer a numerical example to illustrate the construction of quadrature formulas based on scattered data, using ideas described in further detail in \cite{quadconst} in the case of the sphere. 

 Given  a set $\C$ and an integer $n\ge 1$, we define a measure $\nu_\C$ that associates the mass $\lambda_n(\mu;z)$ with each $z\in\C$. If $\{q_k\in\Pi_k\}$ is the (finite) system of orthonormal polynomials with respect to $\nu_\C$, there exist constants $d_{j,k}$ such that $q_k=\sum_j d_{j,k}p_j(\mu)$. If $P\in\Pi_n$ then with $m_0=\sqrt{\mu([-1,1])}$,
\begin{eqnarray*}
\int_{-1}^1 P(x)d\mu(x) &=&\int_{-1}^1\int_{-1}^1 P(t)\sum_{k=0}^n q_k(x)q_k(t)d\nu_\C(t)d\mu(x)=
m_0\int_{-1}^1 P(t)\sum_{k=0}^n d_{0,k}q_k(t)d\nu_\C(t) \\
&=&\sum_{z\in\C}\left\{m_0\lambda_n(\mu;z)\sum_{k=0}^n d_{0,k}q_k(z)\right\}P(z).
\end{eqnarray*}
If $G$ is the matrix with $G_{j,k}=\int_{-1}^1 p_j(\mu;t)p_k(\mu;t)d\nu_\C(t)$, and $D$ is the matrix with $(j,k)$-th entry $d_{j,k}$, then the orthonormality of the system $\{q_k\}$ leads to $G^{-1}=DD^T$. Hence, 
if ${\bf e}=(1,0,\cdots,0)^T\in\RR^{n+1}$ and ${\bf b}=(b_0,\cdots,b_n)$ satsifies $G{\bf b}={\bf e}$, then our quadrature formula is
$$
\int_{-1}^1 P(x)d\mu(x)=\sum_{z\in\C}\left\{m_0\lambda_n(\mu;z)\sum_{j=0}^n b_jp_j(\mu;z)\right\}P(z).
$$
This algorithm will work as long as $G$ is a positive definite matrix. In this case, the system $G{\bf b}={\bf e}$ can be solved using state of the art methods such as the conjugate gradient method. The speed of convergence of this method as well as the quality of the quadrature formula depends upon the condition number of $G$, which in turn, is the square of the ratio of the upper and lower bounds in the M--Z inequalities \eref{polynormequiv} with $\nu_\C$ in place of $\nu$ and $p=2$. While Gram matrices are  ill conditioned in general, they seem to be quite well conditioned for the purpose of this application as long as the degree of the polynomial is not too large in relation with the number of points considered.

As a numerical experiment, we chose 1024 points on the interval $[-1,1]$, chosing one point randomly from each of the 1024 equal subintervals of $[0,\pi]$, and taking the cosine transform to obtain points on $[-1,1]$. For various values of $n$, we attempted to obtain quadrature formulas exact for integrating polynomials in $\Pi_n$ with respect to the Lebesgue measure, $\jsup{\mu}{0}{0}$. The error in the construction was measured by the norm of the difference between the identity matrix of order $\lfloor n/2\rfloor$ and the  matrix of computed inner products $(\sum_{z}w_z\jsup{p_k}{0}{0}(z)\jsup{p_\ell}{0}{0}(z))$, $k,\ell=0,\cdots,\lfloor n/2\rfloor$. Each computation was carried out 30 times using Matlab 7.1 running on Intel Core 2 CPU with 2.13GHz, and 3GB RAM.  The averages of the various parameters of interest are summarized in Table~\ref{quadexpt}. It is interesting to note that we were able to obtain suprizingly accurate quadrature formulas even for degree 1023. Even though there were several negative weights in this case, the total variation of the resulting discrete quadrature measure seems to be surprisingly small.

\begin{table}[h]
$$
\begin{array}{|c|c|c|c|c|c|}
\hline
n & \mbox{cond} & \mbox{pos} & \sum |w_z| & \mbox{error} &\mbox{time}\\
\hline
   256  &  1.5342&          1024&      2.0000&   2e-14 &0.4769\\
\hline
512  &  3.3536&          1023.92 &      2.0000&  9e-14&  1.2441\\
\hline
768 &  25.2965&         1012.28 &   2.0088   &2.4e-13& 2.6159\\
\hline
896& 76.8983&        945.60&    2.1471& 3.8e-13   & 3.7419\\
\hline
1023&   11270.5423&   885.64 &   2.6715&   2.75e-12 &9.5317\\
\hline
 \end{array}
$$
\caption{Averages for 30 trials of: $n=$ degree of exactness, cond is the condition number of the original Gram matrix, pos is the number of positive weights, error is the norm of the difference between the matrices as described in the text, (xe-m means $x*10^{-m}$), time is the time required for each computation, in seconds.}
\label{quadexpt}
\end{table}


\end{document}